\documentclass[11pt,leqno]{article}
\usepackage{amsmath, amsthm, amsfonts, amssymb, color}
\usepackage{mathrsfs}
\newtheorem{thm}{Theorem}[section]

\newtheorem{lem}[thm]{Lemma}
\newtheorem{prp}[thm]{Proposition}

\theoremstyle{definition}

\newcommand{\scr}[1]{\mathscr #1}
\definecolor{wco}{rgb}{0.5,0.2,0.3}

\numberwithin{equation}{section}

\def\R{\mathbb R}

\def\d{\text{\rm{d}}}

\def\beg{\begin} 

\def\Ric{\text{\rm{Ric}}}

\def\Ent{\text{\rm{Ent}}}

\def\Ric{\text{\rm{Ric}}}
\def\cut{\text{\rm{cut}}}

\def\C{\scr C}

\def\Pr{\scr P}

\def\supp{\text{\rm supp}}

\def\veps{\varepsilon}
\def\ra{\rightarrow}
\newcommand{\fin}{\hspace*{\fill}$\Box$}

\title{Characterization of non-constant lower bound of Ricci curvature
via entropy inequality on Wasserstein space}

\author{{\bf Jinghai Shao $^{a)}$ and Bo Wu $^{b)}$}\\
\footnotesize {$^{a)}$ School of Mathematical Sciences, Beijing Normal
University, Beijing 100875, China}\\
\footnotesize {$^{b)}$ School  of Mathematical Sciences, Fudan
University, Shanghai 200433, China}\\
{\footnotesize E-mail: shaojh@bnu.edu.cn, wubo@fudan.edu.cn}}

\date{}

\begin{document}

\maketitle

\begin{abstract} When the Ricci curvature of a Riemannian manifold is not lower bounded by a constant, but lower bounded by a continuous function, we  give a new characterization of this lower bound through the convexity of relative entropy on the probability space over the Riemannian manifold.  Hence, we generalize K.T. Sturm and von Renesse's result (Comm. Pure Appl. Math. 2005) to the case with non-constant lower bound of Ricci curvature.
\end{abstract}

\noindent Keyword: Wasserstein space, Ricci curvature, optimal transport map
\vskip 2cm

\section{Introduction}
In the work \cite{RS}, K.T. Sturm and von Renesse give a characterization of lower bound of Ricci curvature of a Riemannian manifold $M$ through the convexity property of relative entropy on the probability space over $M$. There they considered the setting that Ricci curvature of the manifold is lower bounded by a constant. The goal of this paper is to extend their characterization to the case that the Ricci curvature is  bounded below by a continuous function.
The idea of K.T. Sturm and von Renesse's \cite{RS} has been extended in \cite{LV} \cite{St1} \cite{St2} to metric measure space, and this new definition of lower bound of Ricci curvature owns the stability under the convergence of metric measure space.   We refer the reader to the book \cite{Vil2} for more related works on this widely studied topic.

Let $M$ be a smooth connected Riemannian manifold
of dimension $n$. Fix a point $o\in M$ throughout this work. Denote by $\rho(x, y)$ the Riemannian distance between $x,\,y\in M$, and set $\rho_o(x)=\rho(o,x)$ for $x\in M$.  Let $m(\d x)$ be the Riemannian volume measure on $M$.
Let $\Pr(M)$ be the set of all probability measures on $M$. For $\mu,\,\nu\in \Pr(M)$, the $L^2$-Wasserstein distance between them is defined by
\[ W_2(\mu,\nu)=\inf\Big\{\int_{M\times M}\!\!\rho(x,y)^2\pi(\d x,\d y);\pi\in \C(\mu,\nu)\Big\}^{1/2}\]
where $\C(\mu,\nu)$ stands for the set of all couplings to $\mu$ and $\nu$, that is, the set of all probability measures $\pi$ on $M\times M$ with
$\pi(A\times M)=\mu(A)$ and $\pi(M\times A)=\nu(A)$ for every Borel set $A\subset M$.
$\Pr(M)$ endowed with the metric $W_2(\cdot,\cdot)$ is called a Wasserstein space. The properties of $(\Pr(M), W_2)$ are closely related to the properties of $M$. For example, except the discussion on lower bound of Ricci curvature such as in \cite{RS} and present work, \cite{otto} showed the heat flow  on $M$ can be obtained as the gradient flow on $\Pr(M)$ for the relative entropy. This idea has been extended to deal with other kinds of differential equations in \cite{AGS}.

Let $\Pr_2(M)$ be the set of all probability measures $\mu$ on $M$ which is absolutely continuous with respect to (w.r.t.) volume measure and satisfies $\int_M \rho(o, y)^2 \mu(\d y) < \infty$.
This work strongly depends on the optimal transport map on manifolds. In order to introduce our main result, we recall some basic results on the optimal transport map on Riemannian manifold. By \cite{Mc},  for any pair of absolutely continuous probability measures $\mu_0$ and $\mu_1$ in $\Pr_2(M)$, there exists a unique map $F:M\ra M$ such that $\mu_1=(F)_\ast \mu_0:=\mu_0\circ F^{-1}$. Moreover, there exists a function $\phi$ such that $F$ can be expressed in the form $F(x)=\exp_x(\nabla \phi(x))$. Let $F_t(x)=\exp_x(t\nabla \phi(x))$ for $t\in [0,1]$. Then $\mu_t:=(F_t)_\ast\mu_0$ is the unique geodesic in $(\Pr(M), W_2)$ joining $\mu_0$ to $\mu_1$. By \cite[Proposition 5.4]{CMS}, for each $t\in (0,1)$, $\mu_t$ is absolutely continuous w.r.t. the volume measure $m$ on $M$. Hence, we have $\mu_t\in \Pr_2(M)$, $t\in (0,1)$.

The entropy is defined as a function on $\Pr(M)$ by
$$\Ent(\nu):=\int_M\frac{\d\nu}{\d m}\log\Big(\frac{\d\nu}{\d m}\Big) \d m(x)$$
if $\nu$ is absolutely continuous  w.r.t. volume measure $m$ on $M$  and \[\int_M\max\big\{\log\big(\frac{\d\nu}{\d m}\big), 0\big\}\d \nu(x)<\infty;\]
otherwise,  $ \Ent(\nu):= +\infty$.

Our main result of this work is:
\begin{thm}\label{t1.1} Let $\Ric_x$ be the Ricci curvature at point $x\in M$. Let $x\mapsto K_x$ be a continuous function on $M$. For $D\subset M$, denote by $K(D)=\sup_{z\in D} K_z$. Let $B_x(r)$ denote the open ball of radius $r>0$ centered at $x\in M$. The following
properties are equivalent:
\begin{itemize}
\item[$\mathrm{(i)}$]\  $\Ric_x\geq  -K_x, \quad \forall\,x\in M $.

\item[$\mathrm{(ii)}$]\ For each pair of $\mu_0,\mu_1\in\Pr_2(M)$, it holds: $\forall\, t\in[0,1]$,
\begin{equation}\label{ent-ine}
\begin{split}
\Ent(\mu_t)&\leq (1-t)\Ent(\mu_0)+t\Ent(\mu_1)\\
&\quad +\frac{ t(1-t)}{2}\int_MK\big(B_x( \rho(F(x),x))\big)\rho^2(F (x),x)\mu_0(\d x),
\end{split}
\end{equation}
where the map $F :M\rightarrow M$ is the unique optimal transport map between $\mu_0$ and $\mu_1$ for the $L^2$-Wasserstein distance, which can be expressed by  $F (x) =
\exp_x(-\nabla\varphi)$ for some function $\varphi$.  Here $\mu_t=(F_t)_\ast\mu_0$, $t\in [0,1]$, is the geodesic jointing $\mu_0$ to $\mu_1$, where $F_t(x):=\exp_x(-t\nabla\varphi(x))$ for $t\in [0,1]$.
\end{itemize}
\end{thm}

As an application of this result, we can obtain the following volume growth estimate when the Ricci curvature of Riemannian manifold is locally lower bounded.
\begin{prp}\label{t1.2}
Assume $\Ric_x\geq -K_x$ for every $x\in M$. For fixed $x_0\in M$, consider the volume  $V_R:=m(\bar B_R(x_0))$ of the closed ball centered at $x_0$ with diameter $R$. Then
for all $R\geq 2\veps>0$,
\begin{equation}\label{vol-ine}
V_R\leq V_{2\veps}\Big(\frac{V_{2\veps}}{V_{\veps}}\Big)^{\frac{R}{\veps}}
\exp\Big[\frac{K\big(B_{x_0}( R+2\veps)\big)}{2}R(R+\veps) \Big].
\end{equation}
In particular, if $\Ric_x\geq -C(1+\rho_o(x))$, $\forall\,x\in M$, for some constant $C>0$, then
\begin{equation}\label{vol-ine-1}
V_R\leq V_{2\veps}\Big(\frac{V_{2\veps}}{V_{\veps}}\Big)^{\frac{R }{\veps}} \exp\Big[\frac{C(1+\rho_o(x_0)+R+2\veps)}{2}R(R+\veps)\big) \Big].
\end{equation}
\end{prp}

According to \cite{ATW}, the local curvature bound $\Ric_x\geq -K_x$ can also be characterized by the log-Harnack inequality and gradient estimate of the heat semigroup. To be more precise, let $p_t(x,y)$ be the heat kernel on $M$, that is, the minimal positive functional solution of the heat equation $\big(\frac 12 \Delta-\frac{\partial}{\partial t}\big)p_t(x,y)=0$. Let $(W_t(x))$ be the Brownian motion on $M$ with starting point $x$ and life time $\zeta(x)$, where
$\zeta(x)=\lim_{N\ra \infty}\zeta_N(x)$, and
\[\zeta_N(x):=\inf\{t>0;\rho(W_t(x),o)\geq N\},\ N\geq 1.\]
The associated semigroup is given by
\begin{equation}\label{semi}
P_tf(x)=\mathbb{E}\big[f(W_t(x))\mathbf{1}_{t<\zeta(x)}\big],\quad t\geq 0,
\end{equation}
for $f\in \mathscr B_b(M)$, where $\mathscr B_b(M)$ stands for the set of all bounded measurable functions  on $M$.
For any $D\subset M$, let
\[D_r=\{z\in M;\ \rho(z, D)\leq r\},\quad r>0.\]
For a given bounded open domain $D\subset M$, set
\[\C_D=\{\phi\in C^2(\bar D);\ \phi|_D>0,\ \phi|_{\partial D}=0\}.\]
\begin{prp}[\cite{ATW} Theorem 1.1]
The following statements are equivalent:
\begin{itemize}
  \item[$\mathrm{(i)}$] $\Ric_x\geq -K_x,\quad \forall\, x\in M$.
  \item[$\mathrm{(iii)}$] For any bounded open domain $D\subset M$ and $\phi\in \C_D$, the log-Harnack inequality
      \begin{align*}
        &P_t\log f(y)-\log\Big(P_tf(x)-1-P_t1(x)\Big)\\
        &\leq \frac{\rho(x,y)^2}2\Big(\frac{K(D_{\rho(x,y)})}{1-e^{-2K(D_{\rho(x,y)})t}}
        +\frac{c_D(\phi)^2(e^{2K(D_{\rho(x,y)})t}-1)}{2K(D_{\rho(x,y)})\phi(y)^4}\Big),\\
        &\quad t>0,\ y\in D,\ x\in M,
      \end{align*} holds for strictly positive   $f\in \mathscr B_b(M)$, where
      \[c_D(\phi)=\sup_D\{5|\nabla \phi|^2-
      \frac 12\phi(\Delta \phi)\}\in [0,\infty).\]
  \item[$\mathrm{(iv)}$] For any bounded open domain $D\subset M$ and any $\phi\in \C_D$,
  \begin{align*}
  |\nabla P_t f|^2(x)\leq \{P_t f^2-(P_tf)^2\}(x)\Big(\frac{K(D)}{1-e^{-2K(D)t}}
  +\frac{c_D(\phi)^2(e^{2K(D)t}-1)}{2K(D)\phi(x)^4}\Big)
  \end{align*}holds for all $x\in D$, $t>0$, $f\in \mathscr B_b(M)$.
  If moreover $P_t1=1$, then the statements above are also equivalent to:
  \item[$\mathrm{(v)}$] For any bounded open domain $D\subset M$ and any $\phi\in \C_D$, the Harnack inequality
      \begin{align*}
        &P_tf(y)\leq P_tf(x)\\
        &\quad+\rho(x,y)\Big(\frac{K(D)}{1-e^{-2K(D)t}}
        +\frac{c_D(\phi)^2(e^{2K(D)t}-1)}{2K(D)\inf_{\ell(x,y)}\phi^4}\Big)^{\frac 1 2}\big(P_tf^2(y)\big)^{\frac 12}
      \end{align*}
      holds for non-negative $f\in \mathscr B_b(M)$, $t>0$ and $x,\,y\in D$ such that the minimal geodesic $\ell(x,y)$ linking $x$ with $y$ is contained in $D$.
\end{itemize}
\end{prp}

\section{Proof of the results}
This section is devoted to the proof of our results.
The volume distortion coefficient $v_t$ plays a crucial role in our analysis. We recall the definition and some basic properties of $v_t$ from \cite{CMS}.
For $t\in [0,1]$ and $x,\,y\in M$, let
 \[Z_t(x,y)=\{z\in M; \rho(x,z)=t\rho(x,y)\ \text{and}\ \rho(z,y)=(1-t)\rho(x,y)\}.\]
$Z_t(x,y)$ is the barycenter between $x$ and $y$. For a set $Y\subset M$, define
\[Z_t(x,Y)=\bigcup_{y\in Y} Z_t(x,y).\]
Letting $B_r(y)\subset M$ denote the open ball of radius $r>0$ centered at $y\in M$, for $t\in (0,1]$ the volume distortion coefficient $v_t$ is defined by
\begin{equation}\label{distort}
v_t(x,y)=\lim_{r\ra 0}\frac{m(Z_t(x,B_r(y)))}{m(B_{tr}(y))}.
\end{equation}
It always holds $v_t(x,y)>0$ and $v_1(x,y)=1$.
We recall the following comparison bound on volume distortion from \cite[Corollary 2.2]{CMS}
\begin{lem}[\cite{CMS}]\label{lem-dis}
Assume that $\Ric\geq (n-1)k$ throughout $M$ for some $k\in \R$. Then for $x,y\in M$ with
$y\notin \cut(x)$ and $t\in (0,1)$,
\[v_t(x,y)\geq \left(\frac{S(t\rho(x,y);k)}{S(\rho(x,y);k)}\right)^{n-1},\]
where $S(r;k)$ is defined by
\begin{equation}\label{s-k}
  S(r;k)=\begin{cases}
    \sin\Big(r\sqrt{\frac k{n-1}}\Big)\Big\slash \Big(r\sqrt{\frac k{n-1}}\Big),& k>0,\\
    1, &k=0,\\
    \sinh\Big(r\sqrt{\frac{-k}{n-1}}\Big)\Big\slash \Big(r\sqrt{\frac{-k}{n-1}}\Big), &k<0.
  \end{cases}
\end{equation}
\end{lem}

\noindent\textbf{Proof of Theorem \ref{t1.1}:}\\
 (i)\ $\Rightarrow$\ (ii):
As $\mu_t=(F_t)_\ast \mu_0$, we have for any bounded continuous function $g(\cdot)$ on $M$,
\[\int_{M}g(x)\d \mu_t(x)=\int_M g(F_t(x))\d\mu_0(x),\]
which implies
\[\int_M g(x)\frac{\d \mu_t}{\d m}(x)\d m(x)=\int_M g(F_t(x))\frac{\d \mu_0}{\d m}(x)\d m(x).\]
By changing of variable, we get
\[\int_M g(F_t(x)) \frac{\d \mu_t}{\d m}(F_t(x))\det\big(\d F_t(x)\big)\d m(x)=\int_M g(F_t(x))\frac{\d \mu_0}{\d m}(x)\d m(x).\]
By the arbitrariness of $g$, we get
\[\frac{\d\mu_t}{\d m}(F_t(x))\det\big(\d F_t(x)\big)=\frac{\d \mu_0}{\d m}(x).\]
 Let $J_t(x)=\mathrm{det}(\d F_t(x))$.
It holds
\begin{align*}
  \Ent(\mu_t)&=\int_M  \log \frac{\d \mu_t}{\d m}(x) \d \mu_t(x)\\
  &=\int_M \log\frac{\d\mu_t}{\d m}\big(F_t(x)\big)\d \mu_0(x)\\
  &=\int_M\Big(\log\frac{\d\mu_0}{\d m}(x)-\log J_t(x)\Big)\d \mu_0(x)\\
  &=\Ent(\mu_0)-\int_M \log J_t(x)\d\mu_0(x).
\end{align*}
According to \cite[Theorem 4.2]{CMS}, $F(x)$ does not belong to the cut locus of $x$ for $\mu$ almost everywhere $x$. Then by Lemma 6.1 of \cite{CMS}, 
$J_t=\det( \d F_t(x))$ satisfies the inequality
\begin{equation}\label{ine-J}J_t^{  1/n}(x)\geq (1-t)[v_{1-t}(F(x),x)]^{1/n}+t[v_t(x, F(x))]^{1/n}J_1^{1/n}(x).
\end{equation}
Therefore,
\begin{equation}\label{ine-ent}
\begin{split}
&-\Ent(\mu_t)+(1-t)\Ent(\mu_0)+t\Ent(\mu_1)\\
&=\int_M\log J_t(x)\mu_0(dx)-t\int_M\log J_1(x)\mu_0(dx)\\
&\geq n\!\int_M\!\!\log\big[(1\!-\!t)v_{1-t}(x,F(x))^{\frac 1n}\!+\!tv_t(x,F(x))^{\frac 1n}J_1(x)^{\frac 1n}\big]\mu_0(dx)\\
&\quad -\!t\!\int_M\!\log J_1(x)\mu_0(dx).
\end{split}
\end{equation}
In order to use the comparison bound on volume distortion (\cite[Corollary 2.2]{CMS}), for each $x\in \supp \mu_0$, we consider the open set $B_x(\rho(x,F(x))+\veps)$, which is a geodesic ball centered at $x$ with diameter $\rho(x,F(x))+\veps$, where $0<\veps<1$.   We look on  $B_x( \rho(x,F(x))+\veps)$ as a  manifold, which satisfies
\begin{align*}
 \Ric_z&\geq  -K_z\\
 &\geq  -K\big(B_x(\rho(x,F(x)) +\veps)\big),\quad \forall\ z\in B_x(\rho(x,F(x))+\veps).\end{align*}
Then by noting that the lower bound of Ricci curvature is negative, Lemma \ref{lem-dis}  implies
\begin{equation}\label{e-1}
\begin{split}
&\!\log\big((1\!-\!t)v_{1-t}(x,F(x))^{\frac 1n}\!+\!tv_t(x,F(x))^{\frac 1n}J_1(x)^{\frac 1n}\big)\\
&\!\geq\! \log\bigg((1\!-\!t)\bigg[\frac{S\big((1\!-\!t)\rho(x,F(x));-K(B_x(\rho(x, F(x)+\veps))\big)}{S\big(\rho(x,F(x));-K(B_x(\rho(x,F(x))+\veps)\big)}
\bigg]^{1-\frac  1n }\\
&\quad
 + t\bigg[\frac{S\big(t\rho(x,F(x));-K(B_x(\rho(x, F(x))+\veps)\big)}{S\big(\rho(x,F(x));-K(B_x(\rho(x,F(x))+\veps)\big)}
\bigg]^{1-\frac1n}J_1(x)^{\frac 1n}\bigg).
\end{split}
\end{equation}
Therefore, combining (\ref{e-1}) with (\ref{ine-ent}) and letting $\veps\rightarrow0$, we get
{\small\begin{align*}
&-\Ent(\mu_t)+(1-t)\Ent(\mu_0)+t\Ent(\mu_1)\\
&\geq n\!\int_M\log\bigg[(1-t)\bigg[\frac{S\big((1-t)\rho(x,F(x));-K(B_x(\rho(x,F(x)))\big)}{S\big( \rho(x,F(x));-K(B_x(\rho(x,F(x)))\big)}
\bigg]^{1-\frac1n}\\
&\quad
 + t\bigg[\frac{S\big(t\rho(x,F(x));-K(B_x(\rho(x,F(x)))\big)}{S\big(t\rho(x,F(x)); -K(B_x(\rho(x,F(x)))\big)}
\bigg]^{1-\frac1n}J_1(x)^{\frac1n}\bigg]\mu_0(\d x)\!-\!t\int_M\log J_1(x)\mu_0(\d x)\\
&\geq(n\!-\!1)\!\!\int_M\![(1\!-\!t)\log S\big((1\!-\!t)\rho(x,F(x));-K(B_x(\rho(x,F(x)))\big)\\ &\quad+ t\log S\big(t\rho(x,F(x));-K(B_x(\rho(x,F(x)))\big)
\!-\! S\big(\rho(x,F(x));-K(B_x(\rho(x,F(x)))\big)]\mu_0(dx)\\
&\geq-\frac{ t(1-t)}{2}\int_M K(B_x(\rho(x, F(x)))\rho^2(x,F(x))\mu_0(\d x),
\end{align*}}where we have used the concavity of the logarithm and following inequality (cf. \cite{RS})
$$(1\!-\! t)\log S((1-t)r;k )+t\log S(t r;k)-\log S(r;k)-\frac{t(1-t)}{2} \frac{k}{n-1} r^2\geq 0.$$

(ii)\ $\Rightarrow$\ (i):
If (i) does not  hold, then there exists a point $z_0\in M$ such that $\Ric_{z_0}<-K_{z_0}$. By the continuous property of $K_x$, there exist two positive constants $\varepsilon_0$ and $\delta$ with $\veps_0,\delta<1$  such that
$$\Ric_z<-K_{z_0}-\veps_0,\quad \forall\, z\in B_\delta(z_0),$$
and $B_\delta(z_0)$ is geodesically complete.
Next, similar to \cite{RS}, we can construct two probability measures $\mu_0,\mu_1$ with $\supp(\mu_0)\subset B_\delta(z_0),\supp(\mu_1)\subset B_\delta(z_0)$ such that
\begin{equation}\label{c2}\Ent(\mu_{1/2})-\frac{1}{2}\Ent(\mu_0)-\frac{1}{2}\Ent(\mu_1)\geq \frac{ K_{z_0}+\veps_0/2}{8}W_2^2(\mu_0,\mu_1).\end{equation}
Indeed, let $e_1,e_2,\ldots,e_n$ be an orthonormal basis of $T_{z_0}M$ such that
\[R(e_1,e_i)e_1=k_i e_i \ \ \text{for some numbers $k_i$, $i=1,\ldots,n$}.\]
Then $\sum_{i=1}^nk_i=\Ric_{z_0}(e_1,e_1)\leq -K_{z_0}-\veps_0$.
For $\beta,\,r>0$, let $A_1=B_{\beta}(\exp_{z_0}(re_1))$, $A_0=B_{\beta}(\exp_{z_0}(-r e_1))$ be geodesic balls and \[A_{1/2}=\exp\Big(\big\{y\in T_{z_0} M:\sum_{i=1}^n (y_i/\beta_i)^2\leq 1\big\}\Big)\] with $\beta_i=\beta(1+r^2(k_i+\frac{\veps_0}{2n})/2)$.
By choosing $\beta\ll r\ll \delta$, one gets that $\gamma_{1/2}\in A_{1/2}$ for each minimizing geodesic $\gamma:[0,1]\ra M$ with $\gamma_0\in A_0$, $\gamma_1\in A_1$. Let $\mu_0$ and $\mu_1$ be the normalized uniform distribution in $A_0$ and $A_1$ respectively, and let $\nu$ be the normalized uniform distribution in $A_{1/2}$.
Then
\[\Ent(\mu_0)=\Ent(\mu_1)=-\log m(A_0)=-\log c_n-n\log \beta+O(\beta^2),
\]where $c_n=m(B_1)$ in $\R^n$, and
\begin{align*}
  \Ent(\nu)&=-\log m(A_{1/2})=-\log c_n-\sum_{i=1}^n\log \beta_i+O(\beta^2)\\
  &=-\log c_n-n\log \beta-r^2( \veps_0/2+\sum_{i=1}^n k_i)/2+O(r^4)+O(\beta^2)\\
  &\geq -\log c_n-n\log \beta+\frac{r^2(K_{z_0}+\veps_0/2)}{2} +O(r^4)+O(\beta^2)
\end{align*}
Since the support of $\mu_{1/2}$ is contained in the set $A_{1/2}$, one gets
\[\Ent(\mu_{1/2})\geq \Ent(\nu).\]
Consequently,
\begin{equation} \label{c3}
\begin{split}
  &\Ent(\mu_{1/2})-\frac 12\Ent(\mu_0)-\frac 12 \Ent(\mu_1)\\
  &\geq \frac{r^2(K_{z_0})+\veps_0/2)}{2} +O(r^4)+O(\beta^2)\\
  &\geq \frac{K_{z_0}+\veps_0/2}{8} W_2(\mu_0,\mu_1)^2
\end{split}
\end{equation} for $\beta\ll r\ll \delta$. Hence, we obtain (\ref{c2}).

On the other hand,
since $\supp(\mu_0)$ and $\supp(\mu_1)$ are contained in $B_{r+\beta}(z_0)$, one gets that for each $x\in \supp \mu_0$, $\rho(x,F(x))\leq 2(r+\beta)$. By (ii), we have
\begin{equation}\label{c3}
\begin{split}
&\Ent(\mu_{1/2})-\frac{1}{2}\Ent(\mu_0)-\frac{1}{2}\Ent(\mu_1)\\
&\leq \frac{1}{8}\int_{A_0}K(B_x(\rho(x,F(x)))\rho^2(x,F(x))\mu_0(\d x)\\
&\leq \frac{1}{8}K(B_{z_0}(3 \beta+3r ))W_2^2(\mu_0,\mu_1)\leq \frac 18K(B_{z_0}(6r))W_2^2(\mu_0,\mu_1).
\end{split}
\end{equation}
Since $\lim_{r\ra0} K(B_{z_0}(6r))=K_{z_0}<K_{z_0}+\veps_0/2$,  we can choose $r>0$ small enough so that \[\frac{K(B_{z_0}(6r))}{8}<\frac{K_{z_0}+\veps_0/2}{8}.\]
Thus, by (\ref{c3}) we get
$$\Ent(\mu_{1/2})-\frac{1}{2}\Ent(\mu_0)-\frac{1}{2}\Ent(\mu_1)
<\frac{K_{z_0}+\veps_0/2}{8}W_2^2(\mu_0,\mu_1),$$
which  is in contradiction with (\ref{c2}). We complete the proof of this theorem.
\fin

Next, we shall show how to use our characterization of lower bound of Ricci curvature through convexity of relative entropy to study the volume growth property of Riemannian manifold $M$.

\noindent\textbf{Proof of Proposition \ref{t1.2}:}
Let $\mu_0$ and $\mu_1$ be the uniform distribution on $\bar B_{\veps}(x_0)$ and $\bar B_{R}(x_0)$ respectively, i.e.
\[\mu_0(\d x)=\frac{\mathbf 1_{\bar B_{\veps}(x_0)}}{m(\bar B_{\veps}(x_0))},\quad \mu_1(\d x)=\frac{\mathbf 1_{\bar B_{R}(x_0)}}{m(\bar B_{R}(x_0))}.\]
Let $(\mu_t)_{t\in[0,1]}$ be the geodesic in $\Pr_2(M)$ connecting $\mu_0$ and $\mu_1$. Then, according to Theorem \ref{t1.1},
\begin{align*}
  \Ent(\mu_t)&\leq (1 - t)\Ent(\mu_0) + t\Ent(\mu_1)\\ &\quad + \frac{t(1 - t)}{2}\!\int_M\!K(B_x(\rho(x,F(x)))\rho^2(x,F(x))\mu_0(\d x).
\end{align*}
Since $\mu_t=(F_t)_\ast \mu_0$, we know that $\supp\,\mu_t\subset\bar B_{\veps+t(R+\veps)}(x_0)$. By Jensen's inequality,
\[\Ent(\mu_t)\geq \Ent\Big(\frac{\mathbf 1_{\bar B_{\veps+t(R+\veps)}}}{m(\bar B_{\veps+t(R+\veps)})}\Big)=-\log V_{\veps+t(R+\veps)}.\]
Hence,
\begin{align*}
-\log V_{\veps+t(R+\veps)}&\leq (1-t)\Ent(\mu_0)+t\Ent(\mu_1)\\
&\quad+\frac{ t(1-t)}{2}\int_M\!K(B_x(\rho(x,F(x)))  \rho(x, F(x))^2\d \mu_0(x)\\
&\leq -(1-t)\log V_{\veps}-t\log V_R +\frac{ t(1-t)}{2} K(B_{x_0}( R+2\veps))(R+\veps)^2.
\end{align*}
Take $t=\veps/(R+\veps)$, then  $\veps+t(R+\veps)=2\veps$, and the desired inequality (\ref{vol-ine}) follows immediately from previous inequality.
\fin

\noindent \textbf{Acknowledgments}
This work is supported by NSFC (No.11301030, 11371099), 985-project, Specialized Research Fund for the
Doctoral Program of Higher Education of China (No. 20120071120001).

\beg{thebibliography}{99}
\bibitem{AGS} L. Ambrosio, N. Gigli, G. Savar\'e,   Gradient Flows in Metric Spaces and in the Space of Probability Measures. Lectures in Mathematics ETH Z\"irich. Birkh\"auser Verlag, Basel,  2005.

\bibitem{ATW} M. Arnaudon, A. Thalmaier, F.-Y. Wang, Equivalent Harnack and gradient inequalities for pointwise curvature lower bound.
    Bull. Sci. math., 138 (2014), 643-655.

\bibitem{Chavel} I. Chavel, Riemannian geometry: a modern introduction, Cambridge University Press, 1993.

\bibitem{CMS} Dario Cordero-Erausquin, R. McCann, M. Schmuckenschl\"ager, A Riemannian interpolation inequality \`a la Borel, Brascamb and Lieb. Invent. Math., 146 (2001), 219-257.

\bibitem{otto} R. Jordan, D. Kinderlehrer, F. Otto, The variational formulation of the Fokker-Planck equation. SIAM J. Math. Anal. 29 (1998), 1-17.

\bibitem{Lott} J. Lott, Manifolds with quadratic curvature decay and fast volume growth. Math. Ann. 325 (2003), 525-541.

\bibitem{LV} J. Lott, C. Villani, Ricci curvature for metric-measure spaces via optimal transport, Ann. Math. 169 (3) (2009),  903-991.

\bibitem{Mc} R.J. McCann, Polar factorization of maps on Riemannian manifolds, Geom. Funct. Anal. 11 (2001), 589-608.

\bibitem{RS} K.T. Sturm, von Renesse, Transport inequalities, gradient estimates, entropy, and Ricci curvature. Comm. Pure Appl. Math. 58 (2005), 923-940.

\bibitem{St05} K.T. Sturm, Convex functionals of probability measures and nonlinear diffusions on manifolds, J. Math. Pures Appl. 84 (2005), 149-297.

\bibitem{St1} K.T. Sturm, On the geometry of metric measure spaces. I, Acta Math., 196 (1) (2006), 65-131.

\bibitem{St2} K.T. Sturm, On the geometry of metric measure spaces. II, Acta Math., 196 (1) (2006), 133-177.


\bibitem{Vil2} C. Villani, Optimal Transport, Old and New, Grundlehren Math. Wiss., vol. 338, Springer, Berlin-Heidelberg,  2009.


\end{thebibliography}
\end{document}